\documentclass[12pt]{article}
\usepackage[utf8]{inputenc}
\usepackage[T1]{fontenc}
\usepackage{amsmath, amssymb}
\usepackage{graphicx}
\usepackage{makecell}
\usepackage{xcolor}
\usepackage{tcolorbox}
\usepackage{booktabs}
\usepackage{tabularx}
\usepackage{placeins}
\usepackage[hyphens]{url}
\usepackage[margin=2.5cm]{geometry}

\newcolumntype{C}[1]{>{\centering\arraybackslash}p{#1}}

\title{Compartmental Epidemiological Models: Discussing Fundamental Elements via the SIR Model}

\author{%
Ligia Liani Barz, José Rafael Santos Furlanetto, Fernando Deeke Sasse \\
\normalsize{Departamento de Matemática, Centro de Ciências Tecnológicas} \\
\normalsize{Universidade do Estado de Santa Catarina, Joinville, Brazil} \\
\small{\texttt{{ligia.barz@udesc.br, jose.furlanetto@udesc.br, fernando.sasse@udesc.br}}}
}
\date{}

\begin{document}

\maketitle

\begin{abstract}
We present a didactic treatment of compartmental disease spread modeling centered on the SIR framework, aimed at mathematics education. Our focus is on how modeling choices encode transmission mechanisms and data-collection practices. We begin by stating the structural hypotheses of compartmental models, including population partitioning, homogeneous mixing, flow closure, and Markovian transitions. We then derive the basic SIR dynamic equations in both continuous-time (differential equations) and discrete-time (difference equations) settings. A key goal is to clarify the distinction between density-dependent and frequency-dependent transmission. We do this by constructing the force of infection from the contact rate, prevalence, and per-contact transmissibility while tracking units and interpretations. We go beyond the usual mass balance SIR model by introducing the discrete Kermack–McKendrick perspective through the cumulative force of infection. This naturally leads to discrete SIR updates with an exponential escape probability term. We also derive Bernoulli-based exponential formulations for both density- and frequency-dependent transmission, interpret these probabilistically over observation intervals, and obtain the associated discrete infectivity kernel. This highlights how geometric infectious survival emerges from per-step removal. The paper includes reproducible, introductory-level Python code to implement the basic discrete SIR model, the continuous SIR model, and the exponential discrete variants. We provide instructors with an example of parameter estimation from prevalence data using least squares, along with guidance on interpretation and limitations. This framework is designed to help teachers integrate model assumptions, units, simulation, and inference in a coherent introduction to infectious disease modeling.
\end{abstract}
\vspace{0.5cm}
\section{Introduction}
\label{sec-intro}

Mathematical modeling of epidemics, originally developed for public health, now also has considerable pedagogical value. Educators increasingly face the challenge of explaining terms like transmission rate, prevalence, isolation, and basic reproduction number, which often arise in public discussions but are rarely accompanied by consistent quantitative explanations. Simple compartmental models, especially the SIR (Susceptible–Infectious–Recovered) model, offer teachers an accessible framework for discussing modeling assumptions, parameters, and limitations with learners.

This article is designed to assist educators who wish to: (i) introduce foundational ideas in epidemiological modeling without advanced prerequisites in differential equations; (ii) conceptually and computationally compare continuous and discrete approaches to compartmental models; and (iii) explore relevant aspects of transmission mechanisms, focusing on density- and frequency-dependent types. While introductory literature is vast, many key aspects for effective instruction are addressed only in advanced works; our aim is to close this gap for teachers.

In this work, we emphasize the discrete formulation. We compare two discrete SIR models. The first uses the Euler discretization of differential equations, a linear step approximation. The second uses an exponential transmission term, derived from the Kermack–McKendrick formalism via cumulative force of infection (see \cite{Diekmann2021}). This procedure offers a direct probabilistic interpretation as probability of escaping infection in the interval—and generally avoids negative outputs with a unit time step.

We focus exclusively on epidemics in closed, homogeneous human populations, making the content accessible for classroom instruction. Extensions for non-human groups, as well as demographic, seasonal, age-structured, and stochastic effects, are reserved for future work. To facilitate use by educators, Python code is provided that emphasizes clarity and ease of modification for investigating various models.

Section~\ref{sec_related} gives a literature review of articles in the pedagogical context and lists the main contributions of this work. Section~\ref{sec_compartmental} introduces the fundamental concepts of compartmental epidemiological analysis and explains the structural hypotheses that underpin this approach. This sets the conceptual framework for the article. Section~\ref{sec_SIR_basic} presents the basic SIR model and its dynamical equations occurring in continuous and discrete forms. In Section~\ref{sec_transm}, we use the form of the force of infection to derive two common types of transmission: density-dependent and frequency-dependent. We show the corresponding dynamic equations in both continuous and discrete approaches. Section~\ref{sec_R0} derives expressions for the basic reproduction number for each transmission hypothesis and stresses its role as a threshold parameter. We also discuss its use in the literature. Section~\ref{sec_KM} examines the general Kermack–McKendrick formalism and introduces the cumulative force of infection. This leads to discrete versions with an exponential transmission term. We then derive discrete representations from Bernoulli mechanisms for both density- and frequency-dependent cases, showing that this approach gives a better interpretation of SIR models and provides a basis for generalizations. In Section~\ref{sec_mod_cont_disc}, we systematically compare continuous and discrete approaches, examine the probabilistic interpretation per interval, and discuss computational stability and simulation. Section~\ref{sec_python} presents Python implementations for the various formulations, permitting readers to reproduce simulations and compare results for different parameters and schemes. Section~\ref{sec_estimation} gives practical examples of SIR parameter estimation using real data. Section~\ref{sec_conclusoes} suggests future work.

\section{Related Work and Pedagogical Context}
\label{sec_related}

This work incorporates both pedagogical and non-pedagogical references beyond those listed here.
Our aim is to supplement articles, books, and online materials that teach compartmental epidemiological models, highlighting those described below.

The article by Blackwood et al. \cite{Blackwood2018} provides an introduction to compartmental modeling for university-level students and instructors, along with an overview of model extensions. The discussion on the mathematical form of the force of infection is included, but only frequency-dependent transmission is considered. The book by Vynnycky and White \cite{Vynnycky2010} is often adopted in public health and epidemiology courses, bringing together exercises and case studies; it also discusses the distinction between frequency-dependent and density-dependent transmission in the context of continuous models. Keeling and Rohani \cite{Keeling2008} include, in addition to the continuous approach, a succinct discrete-time formulation. The online resource by Ledder \cite{Ledder2020} organizes classroom activities and non-computational simulation proposals related to SIR and SEIR models. The module by Winkel \cite{Winkel2023} offers an introduction to the SIR model with examples of parameter fitting using spreadsheets. Kang \cite{Kang2024} discusses the SIR model, focusing on exact solutions. McBane \cite{McBane2021} proposes chemical kinetics as an analogy to introduce SIR and simulate the model dynamics; Arnold \cite{Arnold2024} uses concrete physical classroom experiments, employing fluid measurements and probability games to simulate the SIR model, and Hart et al. \cite{Hart2021} explore a lab simulation of physical kinetics in a distance teaching context, with support from spreadsheets and graphical representations. Smith \cite{Smith2022} presents SIR and SIRV models with applications to COVID-19 data; the work also cites a difference equation formulation but does not detail numerical procedures. In \cite{Arjona2022}, the SIR model is presented by analogy with the dynamics of an ethanol oxidation reaction. Meyer et al. \cite{Meyer2023} provide a mathematical introduction to SIR models and their generalizations, covering continuous and discrete formulations and discussing the equilibrium points of more general classes. Finally, Okabe et al. \cite{Okabe2020} derive the SIR model from a growth law and discuss the form of the density-dependent transmission term; the frequency-dependent case is not treated. These authors also present exact solutions and numerical approximations.
A synthesis of this review is presented in Table~\ref{q1}.
\begin{table}[htbp]
\centering
\caption{Pedagogical references.}
\label{q1}
\footnotesize
\setlength{\tabcolsep}{3pt}
\renewcommand{\arraystretch}{1.10}

\begin{tabular}{@{}p{6.1cm}@{\hspace{8pt}}c@{\hspace{10pt}}c@{\hspace{10pt}}c@{}}
\toprule
\textbf{Reference} &
\makecell[c]{\scriptsize\textbf{Differential}\\[-1pt]\scriptsize\textbf{Equations}} &
\makecell[c]{\scriptsize\textbf{Difference}\\[-1pt]\scriptsize\textbf{Equations}} &
\makecell[c]{\scriptsize\textbf{Computational}\\[-1pt]\scriptsize\textbf{Implementation}} \\
\midrule
Arnold et al. (2024) \cite{Arnold2024} & Yes & No & No \\
Blackwood et al. (2018) \cite{Blackwood2018} & Yes & No & No \\
Vynnycky et al. (2010) \cite{Vynnycky2010} & Yes & No & No \\
Keeling et al. (2008) \cite{Keeling2008} & Yes & No & No \\
Ledder (UNL) (2020) \cite{Ledder2020} & No & Yes & No \\
Winkel (QUBES) (2023) \cite{Winkel2023} & Yes & Yes & Yes \\
Kang (2024) \cite{Kang2024} & Yes & No & No \\
McBane (2021) \cite{McBane2021} & Yes & No & No \\
Hart et al. (2021) \cite{Hart2021} & Yes & No & Yes \\
Smith et al. (2022) \cite{Smith2022} & Yes & Yes & No \\
Almanza-Arjona et al. (2022) \cite{Arjona2022} & Yes & No & No \\
Meyer et al. (2023) \cite{Meyer2023} & Yes & Yes & No \\
Okabe et al. (2020) \cite{Okabe2020} & Yes & No & Yes \\
\bottomrule
\end{tabular}

\vspace{1.5mm}
\begin{minipage}{0.95\linewidth}
\footnotesize
\textit{Criterion:} ``Yes'' indicates that the topic is explicitly part of the material's scope.
``Computational Implementation'' refers to examples, routines, or explicit software usage presented as part of the material.
\end{minipage}

\end{table}

In most of the works mentioned above, the system dynamics are described continuously, i.e., as systems of differential equations. Here, in addition to the continuous approach, we also introduce the discrete approach, which uses systems of difference equations. In many cases, the discrete approach is preferable to the continuous one. Some reasons are listed next: (a) epidemiological data are generally provided by day or week; (b) the interpretation of flows between compartments can be understood in terms of probability distributions; (c) the qualitative study of the dynamical system as a \textit{map} (discrete dynamical system) allows for the discussion of stability and thresholds intuitively \cite{Elaydi2005, Allen1994, Wacker2020} and (d) dynamics described by difference equations can be implemented with elementary tools (spreadsheets or introductory programming). Regarding the last item, the formulation of compartmental models in discrete time allows for the introduction of fundamental concepts of epidemiological dynamics to audiences that have not yet had contact (or never will) with the theory of differential equations. Furthermore, even in the context of university education, the analysis of discrete dynamical systems is rarely covered in undergraduate courses, which traditionally privilege continuous formulations. In this sense, discrete modeling not only expands the scope of available tools but also helps to fill a formative gap in understanding system dynamics.

There is a scarcity of materials specifically aimed at higher education mathematics, integrating mathematical rigor, computational accessibility, and relevance to the curriculum. We aim to fill several gaps in the pedagogical literature on compartmental mathematical modeling of epidemics. They are listed below:

\begin{enumerate}
\item Lack of detailed discussions about hypotheses in compartmental models (Section \ref{sec_compartmental}).
\item Absence of a more detailed analysis of the discrete dynamics of models (from Section \ref{sec_SIR_basic} onward).
\item Lack of detailed discussions about transmission mechanisms (Section \ref{sec_transm}).
\item Absence of discussion about the exponential SIR model (Section \ref{sec_KM}).
\item Lack of examples of simple computational implementations (Section \ref{sec_python}).
\item Lack of examples showing how the parameter estimation of a compartmental model can be done using data (Section \ref{sec_estimation}).
\end{enumerate}

\section{Compartmental Models}
\label{sec_compartmental}

\textit{Compartmental modeling} is one of the most widely used approaches to describing the dynamics of infectious diseases. For a description of this and other types of epidemiological models, see, for example, \cite{Vynnycky2010}, p. 20, or \cite{Jacquez1996}. An extensive review of compartmental models, with implementations in the R language, can be found in \cite{Tang2020}.

Most classical compartmental models share a set of basic hypotheses. Understanding such hypotheses is necessary to define the scope and limitations of the models. We will list below some of these hypotheses \cite{Vynnycky2010, Brauer2012}:

\begin{enumerate}
\item \textbf{Well-defined population and counting by compartments.} It is assumed that the population of interest is reasonably delimited and that, at each instant, every individual belongs to exactly one compartment. Thus, compartments form a partition of the population (unless one explicitly includes a compartment ``outside the system'', such as deaths or migration).
\item \textbf{Homogeneous mixing within compartments (homogeneity hypothesis).} We assume that individuals in the same compartment are epidemiologically equivalent for the purpose of the model: they are subject to the same average circumstances and have the same probability (or rate) of undergoing transitions. In practice, this amounts to ignoring heterogeneities in age, spatial structure, contact networks, comorbidities, and behavioral differences. The homogeneity hypothesis is justified when seeking a first-approximation model for large populations.
\item \textbf{Dynamics determined by aggregate rates/flows.} Transitions between compartments are described by flows $\Phi_{i\to j}$ that depend only on the aggregate state $(X_1,\ldots,X_m)$ and model parameters. In other words, it is assumed that the system can be closed by equations in terms of compartmental variables, without tracking individuals.
\item \textbf{Absence of additional memory (Markovianity at the compartment level)}. In the simplest case, state transition rates (or probabilities per step) depend only on the current state of the system, not on the time since an individual entered the compartment.
\item \textbf{Determinism versus stochasticity.} Compartmental models can be stochastic (treating infection and recovery as random events) or deterministic (treating variables as mean values). In this work, we consider only deterministic models, i.e., we do not explicitly account for random fluctuations (demographic noise, environmental noise, underreporting). This choice is appropriate when we want to emphasize mechanisms and qualitative structures (thresholds, equilibria, parameter dependence) and when the population is sufficiently large for averages to be informative.
\item \textbf{Demographic closure (when adopted).} Many models assume a constant total population in the time interval of interest, so that $N=\sum_i X_i$ remains constant. When this hypothesis is not appropriate, we must introduce demographic flows (births, natural mortality, disease-induced mortality, migration) as additional terms in the balance. In this article, when dealing with the basic SIR model, this assumption is justified for relatively fast outbreaks.
\end{enumerate}

In epidemiological compartmental models, the population is partitioned into \textit{compartments} that represent the possible epidemiological states of an individual with respect to the disease. In its simplest form, these states might be ``susceptible'' ($S$), ``infectious'' ($I$), and ``recovered'' ($R$), but several extensions are possible: ``exposed/latent'' ($E$), ``asymptomatic'' ($A$), ``hospitalized'' ($H$), ``vaccinated'' ($V$), among others. The \textit{system state} at an instant $t$ is given by the number (or proportion) of individuals in each compartment, and the model describes how these numbers vary over time as a function of \textit{flows} of transition between compartments.

In general, a compartmental model is specified by \cite{Vynnycky2010, Brauer2012, Diekmann2000, Keeling2008, Anderson1991}:
\begin{enumerate}
\item a list of compartments $X_1, X_2, \ldots, X_m$ (for example, $S,I,R$, denoted by SIR or SEIR, SEIRS, etc.);
\item a set of possible transitions $X_i \to X_j$;
\item a quantitative rule for each flow $\Phi_{i\to j}$, i.e., the rate (or amount per step) at which individuals leave $X_i$ and enter $X_j$.
\end{enumerate}
The dynamic equations result from a \textit{mass balance}\footnote{The term ``mass balance'' is used here in analogy to the principle of conservation of mass used in chemical engineering: the variation of mass within a control volume equals input minus output, plus what is generated minus what is consumed by reactions.} in each compartment: the variation in $X_i$ is ``input minus output''. Depending on the formalism, this balance can be expressed by differential equations (continuous approach) or by difference equations (discrete approach).

\section{Basic SIR Model}
\label{sec_SIR_basic}

The deterministic basic SIR model is the prototype of compartmental models. Its  introduction is commonly attributed to Kermack and McKendrick \cite{Kermack1927} in 1927. It consists of three compartments, one for each of the three types of individuals in the population of interest: susceptible, infectious, and recovered.

\begin{figure}[ht]
\centering
\includegraphics[width=0.6\linewidth]{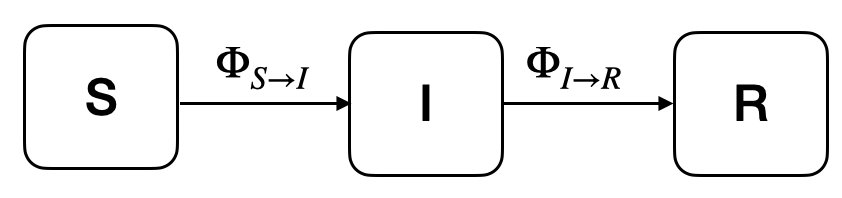}
\caption{Block diagram for the basic SIR model.}
\label{fig_sir1}
\end{figure}

Let us denote, in each compartment, the number of susceptible, infectious, and recovered individuals at a given instant $t$ by $S(t)$, $I(t)$, and $R(t)$, respectively. Such quantities specify the state of the system at a given instant. Furthermore, $N = S(t) + I (t) + R (t)$ denotes the total population, which we will assume to be  \textit{constant} throughout the process under study. In other words, the system is closed. In more general SIR models,  $N$ is not constant.  The flow from compartment $S$ to compartment $I$, i.e., the number of individuals who become infectious per unit time, is denoted by $\Phi_{S\rightarrow I}$. The flow from compartment $I$ to compartment $R$ is denoted by $\Phi_{I\rightarrow R}$. In Figure \ref{fig_sir1}, we show the block diagram of these transitions.

\FloatBarrier

The system dynamics are specified by the flow relationships described in Figure \ref{fig_sir1}. These are generally specified in terms of the system of differential equations \cite{Keeling2008}:

\begin{align}
\label{eq-continua1}
\frac{\mathrm{d}S}{\mathrm{d}t} &= -\Phi_{S\rightarrow I}, \\
\label{continua2}
\frac{\mathrm{d}I}{\mathrm{d}t} &= \Phi_{S\rightarrow I}-\Phi_{I\rightarrow R}, \\
\label{eq-continua3}
\frac{\mathrm{d}R}{\mathrm{d}t} &=\Phi_{I\rightarrow R},
\end{align}
with initial conditions
\begin{equation}
\label{eq-ed4}
S(t_0)=S_0,\quad I(t_0)=I_0,\quad R(t_0)=R_0.
\end{equation}
From now on, we will refer to the approach using differential equations as the \textit{continuous approach}.

Another possible dynamic approach is what we will call \textit{discrete}, expressed in terms of difference equations \cite{Elaydi2005, Martcheva2015}. The discretization of system 
\eqref{eq-continua1}-\eqref{eq-continua3}, with the same initial conditions, using the Euler method, is given by

\begin{align}
\label{eq_discreta1}
S_{n+1}&= S_n-\Phi_{S\rightarrow I} \Delta t, \\
I_{n+1}&= I_n+\Phi_{S\rightarrow I} \Delta t- \Phi_{I\rightarrow R} \Delta t, \\
R_{n+1}&= R_n+\Phi_{I\rightarrow R} \Delta t.
\end{align}
In this formalism, it is usual to take the same time unit (days, hours, etc.) for the time step for which the model parameters are empirically estimated, i.e., $\Delta t =1$, so that the model dynamics are given by

\begin{align}
\label{eq_discreta2a}
S_{n+1}&= S_n-\Phi_{S\rightarrow I}, \\
\label{eq_discreta2b}
I_{n+1}&= I_n+\Phi_{S\rightarrow I} - \Phi_{I\rightarrow R}, \\
\label{eq_discreta2c}
R_{n+1}&= R_n+\Phi_{I\rightarrow R}.
\end{align}

We will discuss the advantages, disadvantages, and consequences of continuous and discrete modeling. At this point, it is interesting to emphasize that the discrete approach can be introduced to students and professionals without knowledge of differential equations, since the construction of the discrete model equations \eqref{eq_discreta2a}–\eqref{eq_discreta2c} can be justified based on the notion of input and output flows between compartments.

The next step in model construction is to specify the flows $\Phi_{S\rightarrow I}$ and $\Phi_{I\rightarrow R}$.

The simplest hypothesis is to assume that the infection flow $\Phi_{S\rightarrow I}$, i.e., the number of individuals infected per unit time, is proportional to the number of susceptible individuals at each instant \cite{Coburn2009}, i.e.,
\begin{equation}
\label{eq-mbe1}
\Phi_{S\rightarrow I}= \lambda S.
\end{equation}
The proportionality factor $\lambda$ (not necessarily constant) is called the \textit{force of infection}. Although the assumption of linearity in $S$, expressed in \eqref{eq-mbe1}, can be understood as a mathematical simplification of the model, it can be discussed, justified, and questioned. One possibility is to interpret \eqref{eq-mbe1} as resulting from the assumption that each susceptible individual becomes infectious according to a Poisson process \cite{Montgomery2011} with parameter $\lambda$. In other words, the time until infection follows an exponential distribution with rate $\lambda$. The probability of a susceptible being infected in a sufficiently small time interval $\Delta t$ is approximately $\lambda \Delta t$. The factor $\lambda S$ is the expected value of new infections per unit time (p. 18 of \cite{Keeling2008}). For models with nonlinear incidence in $S$, see, for example, \cite{Liu1987} and \cite{Hethcote1991}.

On the other hand, we will assume that the recovery flow, denoted by $\Phi_{I\rightarrow R}$, is proportional to the number of infectious individuals, i.e.,

\begin{equation}
\label{eq-mbe7}
\Phi_{I\rightarrow R}=\gamma I,
\end{equation}
where $\gamma$ is the recovery rate, normally considered constant. We can understand $1/\gamma$ as the average number of time units that an infectious individual takes to recover. Another way to understand this hypothesis is to assume that the removal of infectious individuals follows a Poisson process with a constant rate $\gamma$ (probability per unit time). Therefore, the time $T$ that an individual remains in compartment $I$ follows an exponential distribution

\begin{equation}
\label{eq-exp1}
f(t)=\gamma \mathrm{e}^{-\gamma t},
\end{equation}
so that the probability of an individual remaining in $I$ for a time $t$ is

\begin{equation}
\label{eq-exp2}
P(t<T)=\mathrm{e}^{-\gamma t},
\end{equation}
The expected value of $T$ with such a distribution is, in fact, $E(T)=1/\gamma$. In a small time interval $\Delta t$, the probability that an infectious individual recovers is given by

\begin{equation}
\label{eq-exp3b}
P_{rec}=1-\mathrm{e}^{-\gamma \Delta t}\approx \gamma \Delta t.
\end{equation}

If we have $I(t)$ infectious individuals, the expected number of individuals that recover $\Delta R$ in a small  time interval $\Delta t$ is $I(t) P_{rec}$, i.e.,
\begin{equation}
\label{eq-exp3c}
\Delta R = I(t) \gamma \Delta t,
\end{equation}
which leads us to \eqref{eq-mbe7} (p. 15 of \cite{Diekmann2000}).

The model dynamics equations now have, in the continuous approach, the following form:
\begin{align}
\label{eq-continua1b}
\frac{\mathrm{d}S}{\mathrm{d}t} &= -\lambda S, \\
\label{eq-continua2b}
\frac{\mathrm{d}I}{\mathrm{d}t} &= \lambda S-\gamma I, \\
\label{eq-continua3b}
\frac{\mathrm{d}R}{\mathrm{d}t} &=\gamma I.
\end{align}
In the discrete approach,

\begin{align}
\label{eq_discreta3a}
S_{n+1}&= S_n-\lambda S_n, \\
\label{eq_discreta3b}
I_{n+1}&= I_n+\lambda S_n - \gamma I_n, \\
\label{eq_discreta3c}
R_{n+1}&= R_n+\gamma I_n.
\end{align}

The form of the force of infection $\lambda$, given by \eqref{eq-mbe1}, is crucial for the type of model, as it specifies the transmission mechanism. In the next section, we will consider two basic types of transmission.

\section{Basic Types of Transmission}
\label{sec_transm}

Let us consider some common mathematical forms of the force of infection $\lambda$. Although we are in the context of the basic SIR model, the discussions presented here hold for the other SIR model extensions we will see later. Let us assume that the force of infection has the form
\begin{equation}
\label{eq_beta1}
\lambda = cp\nu,
\end{equation}
where $c$ is the contact rate, $p$ is the probability that the contact is, in fact, with an infectious individual \footnote{In this article, the notation $p$ will have a meaning determined by the context.}, and $\nu$ is the probability that the contact event between an infectious individual and a susceptible person leads to transmission.

A simplifying hypothesis is that individual contact events are equally likely within the compartment, i.e., $p = I/N$. Such a hypothesis does not account for social, seasonal, and spatial factors that contribute to population heterogeneity within the compartment. This term is called \textit{infection prevalence}. We then have
\begin{equation}
\label{eq-beta2}
\lambda = c\,\frac{I}{N}\,\nu,
\end{equation}
so that equation \eqref{eq-mbe1}  becomes:
\begin{equation}
\label{eq-beta3}
\Phi_{S\rightarrow I}= c\,\frac{I}{N}\,\nu S.
\end{equation}
This relationship is of the \textit{mass action law} type, analogous to that from chemical kinetics. A discussion of this analogy and its limitations can be found in \cite{Heesterbeek2005}. Even assuming homogeneity, the hypotheses specified by \eqref{eq-mbe1} and \eqref{eq_beta1} imply that $\lambda(t)$ does not change considerably in the time interval $\Delta t$ of infection and that infections occur through a Poisson process. A detailed explanation can be found in \cite{Keeling2008}, p. 18.

The form of the contact rate $c$ defines the two basic types of transmission:

\paragraph{(i) Density-dependent transmission}

This is the type of transmission originally used in the SIR model. Denoting by $A$ the area occupied by the population, we assume that the contact rate is proportional to the population density $N/A$, i.e.,
\begin{equation}
\label{eq-beta4}
c = \kappa \frac{N}{A},
\end{equation}
with $\kappa$ being the proportionality constant. We then have
\begin{equation}
\label{eq-beta4b}
\lambda = \kappa\frac{N}{A}\,\frac{I}{N}\,\nu = \kappa \nu\frac{I}{A}.
\end{equation}
and
\begin{equation}
\label{eq-beta5}
\Phi_{S\rightarrow I}= \kappa\frac{N}{A}\,\frac{I}{N}\,\nu S = \kappa \nu\frac{I S}{A}.
\end{equation}
Defining a new parameter called \textit{the transmission coefficient}, or \textit{infectiousness},

\begin{equation}
\label{eq-beta6}
\beta := \kappa \nu,
\end{equation}
we obtain
\begin{equation}
\label{eq-beta7}
\Phi_{S\rightarrow I}= \frac{\beta}{A} S I.
\end{equation}
Observing the first term of equation \eqref{eq-beta4b}, we note that the force of infection increases with population density $N/A$ and also with the proportion of infectious individuals $I/N$, which justifies the designation \textit{density-dependent transmission}. The resulting flow \eqref{eq-beta7} is appropriate for epidemiological models involving densely populated, spatially limited environments (such as urban areas, schools, etc.) where infections are transmitted by simple physical proximity, airborne transmission, or, more generally, by contacts that increase with population density. This is the case for various types of influenza, tuberculosis, meningococcal meningitis, and measles. Such a model is also appropriate for vector-borne diseases, such as dengue, malaria, and yellow fever, where mosquito density is proportional to human population density.

In the special case where:
\begin{enumerate}
\item The total area $A$ is constant—not only for a given population at a given time, but also between different populations that are compared.
\item The probability $\nu$ that contact between an infectious individual and a susceptible person actually leads to transmission is constant, so $\beta =\kappa \nu$ is constant.
\end{enumerate}
we can define a new transmission coefficient $\bar{\beta}:=\beta/A$ so that
\begin{equation}
\label{eq-beta7b}
\Phi_{S\rightarrow I}= \bar{\beta} S I.
\end{equation}
Normally, $\beta$ (or $\bar{\beta}$) is empirically determined for each type of disease.

Although the form \eqref{eq-beta7} is widely used in the literature to model epidemics, both in human and general populations, it is important to note, as Begon et al. \cite{Begon2002} pointed out, that it is interesting to keep in mind the original equations \eqref{eq-beta5} and \eqref{eq-beta7} in order to avoid inappropriate simplifications. With these considerations, the dynamic equations of the model in the case of density-dependent transmission are given, in the continuous approach, by
\begin{align}
\label{eq-continua_densidade1}
\frac{\mathrm{d}S}{\mathrm{d}t} &= -\overline{\beta} I S, \\
\label{eq-continua_densidade2}
\frac{\mathrm{d}I}{\mathrm{d}t} &= \overline{\beta} I S-\gamma I, \\
\label{eq-continua_densidade3}
\frac{\mathrm{d}R}{\mathrm{d}t} &=\gamma I.
\end{align}
In the discrete approach, we have
\begin{align}
\label{eq_discreta4a}
S_{n+1}&= S_n-\overline{\beta}S_n I_n, \\
\label{eq_discreta4b}
I_{n+1}&= I_n+\overline{\beta}S_n I_n - \gamma I_n, \\
\label{eq_discreta4c}
R_{n+1}&= R_n+\gamma I_n.
\end{align}

\paragraph{(ii) Frequency-dependent transmission}

Another type of transmission used in compartmental models is applicable to cases where the contact rate $c$ equals a constant $\eta$ (independent of population density):

\begin{equation}
\label{eq-beta8}
c = \eta = \text{const.},
\end{equation}
which results in a force of infection, still assuming $p=I/N$,
\begin{equation}
\label{eq-beta8b}
\lambda = \eta p \nu = \eta \frac{I}{N} \nu.
\end{equation}
Equation \eqref{eq-beta3} becomes
\begin{equation}
\label{eq-beta9}
\Phi_{S\rightarrow I}= \eta\frac{I}{N}\nu S.
\end{equation}
Defining a new transmission coefficient (with different units from $\beta$)
\begin{equation}
\label{eq-beta10}
\beta' = \eta \nu,
\end{equation}
we have
\begin{equation}
\label{eq-beta11}
\Phi_{S\rightarrow I}= \frac{\beta' S I}{N}.
\end{equation}
Note that the force of infection $\lambda = \eta \nu I/N$ is proportional to infection prevalence $I/N$ and to the constant contact rate $\eta$, which is referred to as \textit{contact frequency}, justifying the designation \textit{frequency-dependent transmission}. Such a model is appropriate for cases of diseases whose transmission depends on contact associated with behavioral and social factors. The risk of infection depends more on the proportion of infectious individuals in these groups than on the absolute number in the general population. In this context, the number of contacts per person does not grow indefinitely as population density increases, as interactions are limited by social networks and routines. In this category are sexually transmitted diseases (individuals have a limited number of sexual partners, regardless of population density).

The epidemic dynamics in the continuous approach are now given by:
\begin{align}
\label{eq-continua_freq1}
\frac{\mathrm{d}S}{\mathrm{d}t} &= -\frac{\beta' I S}{N}, \\
\label{eq-continua_freq2}
\frac{\mathrm{d}I}{\mathrm{d}t} &= \frac{\beta' I S}{N}-\gamma I, \\
\label{eq-continua_freq3}
\frac{\mathrm{d}R}{\mathrm{d}t} &=\gamma I.
\end{align}
In the discrete approach, we have:
\begin{align}
\label{eq_discreta5a}
S_{n+1}&= S_n-\frac{\beta' S_n I_n}{N}, \\
\label{eq_discreta5b}
I_{n+1}&= I_n+\frac{\beta' S_n I_n}{N} - \gamma I_n, \\
\label{eq_discreta5c}
R_{n+1}&= R_n+\gamma I_n.
\end{align}
The discussion about the choice between density-dependent and frequency-dependent forces of infection in compartmental models is made in several articles and books (see, for example, \cite{Begon2002, Ryder2007, Lockhart1996, McCallum2001}), but the adoption of one of these approaches is often done implicitly in the literature in this area. A model combining density-dependent and frequency-dependent forces of infection was presented in \cite{Ryder2007}. Certainly, for a population with fixed $N$, modeling based on density-dependent or frequency-dependent transmission is mathematically equivalent, resulting only in a change of scale for $\beta$, aside from its meaning. We designate any model with dynamics given by the differential or difference equations discussed above as the \textit{basic SIR model}.
\section{Basic reproduction number}
\label{sec_R0}

 In the case of density-dependent transmission, for an outbreak to begin, the net flow in compartment $I$ must initially be positive, i.e.,
\begin{equation}
\label{eq-mbe8}
\Phi_{S\rightarrow I} - \Phi_{I\rightarrow R}=
\overline{\beta} S I -\gamma I>0.
\end{equation}
or
\begin{equation}
\label{eq-mbe9}
\overline{\beta} S (t_0 ) >\gamma.
\end{equation}
Since at the initial instant $t_0$ we must have $S\approx N$, the condition for an outbreak to start becomes

\begin{equation}
\label{eq-mbe10}
\overline{\mathcal{R}}_0:=\frac{\overline{\beta} N}{\gamma}>1.
\end{equation}
The factor $\overline{\mathcal{R}}_0$ is called the \textit{basic reproduction number}. It can be understood as \textit{the expected number of secondary infections produced by a single infectious individual in a totally susceptible population} \cite{Li2011}. The term \textit{secondary infections} refers to the new infections resulting from a single initial case in a population. In the present case, the number of secondary infections increases with population size $N$. If $\overline{\mathcal{R}}_0<1$, the outbreak does not start. We should note that, as $\overline{\beta}N$ and $\gamma$ represent rates (units of $\mathrm{time}^{-1}$), the quantity $\overline{\mathcal{R}}_0$ is dimensionless.

In the case of the SIR model with frequency-dependent transmission, the term $S\approx N$ in $\Phi_{S\rightarrow I} - \Phi_{I\rightarrow R}$ cancels with the denominator $N$, and we have
\begin{equation}
\mathcal{R}'_0 := \frac{\beta'}{\gamma}.
\end{equation}
The values of the basic reproduction number found in much of the literature suggest that $\overline{\mathcal{R}}_0$ and $\mathcal{R}'_0$ measure the infectiousness of a disease. However, as observed in \cite{Li2011, Delamater2019, Lim2020}, such an interpretation is inadequate, as the value of the basic reproduction number for a given epidemic is highly influenced not only by the model used but also by the assumptions and estimation methods adopted. As this topic is controversial, a comprehensive discussion is outside the scope of this work and will be left for future research.

\section{General Kermack–McKendrick SIR Model}
\label{sec_KM}

\subsection{SIR Model with cumulative force of infection}
\label{subsec_sir_cum}

The basic SIR model presented previously is often referred to as the \textit{Kermack}–\linebreak \textit{McKendrick model}, as it appeared in the seminal article \cite{Kermack1927}. However, as observed in \cite{Breda2012}, Kermack and McKendrick actually  proposed a much more general model. The discussion of some of its aspects is worthwhile since it illuminates hypotheses that may not be explicit in an introductory study.  The basic discrete SIR model, described by \eqref{eq_discreta4a}–\eqref{eq_discreta4c} or \eqref{eq_discreta5a}–\eqref{eq_discreta5c}, has the problem that, in the discrete approach, when we have relatively large values of $\overline{\beta}$ (density-dependent transmission) or $\beta'/N$ (frequency-dependent transmission), we can obtain negative values of $S$. As shown in \cite{Diekmann2021}, this is due to the fact that the relation $S_{n+1}=(1-\lambda) S_n$ does not take into account that an individual can become infected only once. To illustrate the nature of the approximations we use in this work, let us consider, in a general way, the model of an epidemic in a closed population, as originally proposed by Kermack and McKendrick. The fundamental assumption of this model is that the number of susceptible individuals $S$ in the population satisfies the relation
\begin{equation}
\label{eq-km1}
\frac{d S(t)}{dt}=-\lambda(t) S(t),
\end{equation}
where $\lambda(t)$ is the force of infection at time $t$, interpreted as the probability per unit time (instantaneous infection rate) that a susceptible individual becomes infectious at time $t$. Formally integrating \eqref{eq-km1}, we can write
\begin{equation}
\label{eq-km2}
S(t+1) = \mathrm{e}^{-\Lambda(t)} S(t),
\end{equation}
where
\begin{equation}
\label{eq-km3}
\Lambda(t) = \int_t^{t+1}\lambda(\tau) \,d\tau
\end{equation}
is the so-called \textit{cumulative force of infection} on a typical susceptible individual in the interval $(t, t+1]$. The quantity $\mathrm{e}^{-\Lambda(t)}$ is interpreted as the probability that a susceptible individual escapes infection in the interval $(t, t+1]$.

It is important to note that only for small values of $\Lambda$ can we use the approximation
\begin{equation}
\label{eq-km4}
S(t+1) \approx [1-\Lambda(t)]\,S(t),
\end{equation}
which is similar in form to equation \eqref{eq_discreta3a}.

One of the most important points of the fundamental paper by Kermack and McKendrick consists of the observation that the force of infection $\lambda(t)$ at time $t$ must be influenced by all individuals who have been infected at time $\tau<t$. This fact can be expressed in terms of the survival and renewal equation
\begin{equation}
\label{eq-km5}
\lambda(t) = \int_0^{\infty} A(\tau) \lambda(t-\tau)S(t-\tau) \,d\tau,
\end{equation}
where the so-called \textit{infectivity kernel} $A(\tau)$ is the expected contribution to the force of infection of an individual who was infected $\tau$ time units ago \cite{Breda2012}. It is easy to show that, substituting \eqref{eq-km5} into \eqref{eq-km3} and using \eqref{eq-km1}, we obtain
\begin{equation}
\label{eq-km6}
\Lambda(t) = \int_t^{t+1} \lambda(\tau)\, d\tau=
\int_0^{\infty} A(\tau)\Big[S(t-\tau)-S(t+1-\tau)\Big]\,d\tau.
\end{equation}
The discrete version of this equation is given by
\begin{equation}
\label{eq-km7}
\Lambda(t)=\sum_{j=1}^{\infty}A_j\Big[S(t-j)-S(t+1-j)\Big]
=\sum_{j=1}^{\infty}A_j\,C(t-j),
\end{equation}
with $C(t-j)=S(t-j)-S(t+1-j)$ being the incidence in the interval $(t-j,t+1-j]$, i.e., the number of new infections occurring $j$ steps ago.

Relation \eqref{eq-km7} shows that $\Lambda(t)$ is decomposed as the sum of contributions coming from groups of previously infectious individuals. The term $A_j$ can be interpreted as an \textit{infectivity weight by infection age}: it quantifies the expected contribution in the current interval $(t,t+1]$ of an individual who was infected in the interval $(t-j,t+1-j]$. In epidemiological terms,   $\{A_j\}_{j\ge1}$  is a discrete profile of ``infectivity over time since infection''.

In this work, we will not investigate techniques for the empirical determination of $A_j$. However, we will theoretically construct SIR models that, although they are particular cases, are more general than the simple SIR models that use the approximation \eqref{eq-km4}. We will next present the form of the corresponding series ${A_j}$.

\subsection{Density-dependent Bernoulli transmission mechanism}

In the discrete model, suppose the infection mechanism is a Bernoulli process \cite{Montgomery2011}, i.e., during the interval $(t, t+1]$, a single infectious individual infects a given susceptible individual with a probability of $p$. Therefore, the probability that a susceptible individual escapes infection caused by that infectious individual is $1-p$. If we have $I(t)$ infectious individuals at time $t$ and assume that the infection events happen occur independently for the susceptible individual in the interval $(t, t+1]$, then the probability that the susceptible individual escapes infection from all infectious individuals is given by $(1-p)^{I(t)}$. We can rewrite this term as follows:
\begin{equation}
\label{eq-km8}
(1-p)^{I(t)} = \mathrm{e}^{I(t)\ln(1-p)}=\mathrm{e}^{-\overline{\beta} I(t)},
\end{equation}
where
\begin{equation}
\label{eq-km9}
\overline{\beta}:=-\ln(1-p)\geq 0.
\end{equation}
Compared with \eqref{eq-km2}, we have
\begin{equation}
\label{eq-km10}
\Lambda(t) = \overline{\beta} I(t),
\end{equation}
so that
\begin{equation}
\label{eq-km11}
S(t+1) = \mathrm{e}^{-\overline{\beta} I(t)} S(t).
\end{equation}
To complete the SIR model, suppose that the probability that an infectious individual remains infectious until the next time step is $1-\gamma$, so that $\gamma$ is the probability that they are removed to compartment $R$, acquiring permanent immunity. Furthermore, since $\mathrm{e}^{-\overline{\beta} I(t)}$ is the probability that a susceptible individual escapes infection, we must have:
\begin{align}
\label{eq-km12}
S(t+1) &= \mathrm{e}^{-\overline{\beta} I(t)} S(t), \\
\label{eq-km13}
I(t+1) &=(1-\mathrm{e}^{-\overline{\beta} I(t)}) S(t)+(1-\gamma)I(t), \\
\label{eq-km14}
R(t+1)&= \gamma I(t)+R(t).
\end{align}
The corresponding difference equations are
\begin{align}
\label{eq-km15}
S_{i+1} &= \mathrm{e}^{-\overline{\beta}I_i} S_i, \\
\label{eq-km16}
I_{i+1} &=(1-\mathrm{e}^{-\overline{\beta} I_i}) S_i+(1-\gamma)I_i, \\
\label{eq-km17}
R_{i+1}&= \gamma I_i+R_i.
\end{align}
From a computational point of view, it is convenient to define $I_{i+1}$ in terms of the incidence $S_i-S_{i+1}$, so that
\begin{align}
\label{eq-km15b}
S_{i+1} &= \mathrm{e}^{-\overline{\beta} I_i} S_i, \\
\label{eq-km16b}
I_{i+1} &= S_i-S_{i+1}+(1-\gamma)I_i, \\
\label{eq-km17b}
R_{i+1}&= \gamma I_i+R_i.
\end{align}

\subsection{Frequency-dependent Bernoulli transmission mechanism}

To construct the corresponding frequency-dependent transmission model, we must assume that each infectious individual makes a number of potentially infectious contacts per time interval that is approximately independent of the population size $N$. Consequently, for a fixed susceptible individual, the probability $p$ of being infected by a given infectious individual in the interval $(t,t+1]$ must be proportional to $1/N$, i.e.,
\begin{equation}
\label{eq-km18}
p=\frac{k}{N},
\end{equation}
where $0\leq k/N\leq 1$,  $k$ is independent of $N$ and represents the expected number of contacts per step. Therefore, the probability that a susceptible individual escapes infection caused by that infectious individual is $1-k/N$. If we have $I(t)$ infectious individuals at time $t$ and assume that the corresponding infection events occur independently in the interval $(t, t+1]$, then the probability that the susceptible individual escapes infection caused by all infectious individuals is given by $(1-k/N)^{I(t)}$. Therefore,
\begin{equation}
\label{eq-km19}
\left(1-\frac{k}{N}\right)^{I(t)} = \mathrm{e}^{I(t)\ln(1-k/N)}=\mathrm{e}^{-\beta' I(t)/N},
\end{equation}
where
\begin{equation}
\label{eq-km20}
\beta':=-N\ln\left(1-\frac{k}{N}\right)\geq 0.
\end{equation}
Comparing with \eqref{eq-km2} we have
\begin{equation}
\label{eq-km21}
\Lambda(t) = \frac{\beta'}{N} I(t),
\end{equation}
so that
\begin{equation}
\label{eq-km22}
S(t+1) = \mathrm{e}^{-\beta' I(t)/N} S(t).
\end{equation}
The corresponding difference equations are now
\begin{align}
\label{eq-km23}
S_{i+1} &= \mathrm{e}^{-\beta' I_i/N} S_i, \\
\label{eq-km24}
I_{i+1} &=(1-\mathrm{e}^{-\beta'I_i/N}) S_i+(1-\gamma)I_i, \\
\label{eq-km25}
R_{i+1}&= \gamma I_i+R_i.
\end{align}
or, in terms of incidence,
\begin{align}
\label{eq-km23b}
S_{i+1} &= \mathrm{e}^{-\beta' I_i/N}S_i, \\
\label{eq-km24b}
I_{i+1} &= S_i-S_{i+1}+(1-\gamma)I_i, \\
\label{eq-km25b}
R_{i+1}&= \gamma I_i+R_i.
\end{align}
It is worth justifying the choice of definition \eqref{eq-km20}. Mathematically, we could have defined $\beta_N:=-\ln(1-k/N)$, such that $\Lambda(t) = \beta_N I(t)$. However, in that case,
\begin{equation}
\label{eq-km26}
\beta_N = -\ln \left(1-\frac{k}{N}\right)= \frac{k}{N}+O\left(\frac{1}{N^2}\right),
\end{equation}
so that $\beta_N\to 0$ when $N \to \infty$. This would imply the undesirable property of a model parameter being strongly dependent on $N$. On the other hand,
\begin{equation}
\label{eq-km27}
\beta' = -N\ln \left(1-\frac{k}{N}\right)= k+O\left(\frac{1}{N}\right),
\end{equation}
so that $\beta'\to k$ when $N \to \infty$. That is, $\beta'$ is approximately independent of $N$ when $N$ is large and, therefore, can be used in comparisons between different populations.

\subsection{Explicit form of the infectivity weight in the exponential SIR model}
\label{subsec:Ak_exponencial}

In this subsection, we derive the form of the coefficients $A_j$  (discrete infectivity kernel) associated with the discrete SIR model with an exponential term. In the exponential SIR model, which assumes Bernoulli removal at each time step, an infectious individual remains infectious until the next step with probability  $1-\gamma$, and is removed with probability $\gamma$ (cf. \eqref{eq-km14}–\eqref{eq-km17b}). Thus, the dynamics of $I$ can be expressed in terms of $C(t)$ as
\begin{equation}
\label{eq_kmA3}
I(t+1)=C(t)+(1-\gamma)I(t).
\end{equation}
Substituting $t\to t-1$ in \eqref{eq_kmA3}, we obtain
\begin{equation}
\label{eq_kmA5}
I(t)=C(t-1)+(1-\gamma)I(t-1).
\end{equation}
To \textit{iterate backwards} identity \eqref{eq_kmA5}, we substitute the terms $I(t-1), I(t-2), \ldots$, recursively.
From \eqref{eq_kmA5} we have
\begin{equation}
\label{eq_kmA6}
I(t-1)=C(t-2)+(1-\gamma)I(t-2).
\end{equation}
Substituting this expression into \eqref{eq_kmA5} we obtain
\begin{equation}
\label{eq_kmA7}
I(t)=C(t-1)+(1-\gamma)C(t-2)+(1-\gamma)^2 I(t-2).
\end{equation}
On the other hand,
\begin{equation}
\label{eq_kmA8}
I(t-2)=C(t-3)+(1-\gamma)I(t-3),
\end{equation}
so that
\begin{equation}
\label{eq_kmA9}
I(t)=C(t-1)+(1-\gamma)C(t-2)+(1-\gamma)^2 C(t-3)+(1-\gamma)^3 I(t-3).
\end{equation}
Proceeding in this way, after $n$ iterations, we obtain
\begin{equation}
\label{eq_kmA10}
I(t)=\sum_{j=1}^{n}(1-\gamma)^{j-1}\,C(t-j)\;+\;(1-\gamma)^{n}\,I(t-n).
\end{equation}

We assume that, in the past, there exists $t_0$ such that $I(t-n)=0$ whenever $t-n<t_0$. Therefore, taking $n$ sufficiently large (with $t-n<t_0$), the residual term $(1-\gamma)^n I(t-n)$ vanishes, so that
\begin{equation}
\label{eq_kmA10b}
I(t)=\sum_{j=1}^{\infty}(1-\gamma)^{j-1}\,C(t-j).
\end{equation}
As in the exponential SIR model, the update of $S$ has the form (cf. \eqref{eq-km2})
\begin{equation}
\label{eq_kmA12}
S(t+1)=\mathrm{e}^{-\Lambda(t)}S(t)\,
\end{equation}
the choice of transmission mechanism amounts to specifying $\Lambda(t)$ as a function of $I(t)$. Let us treat the cases of density-dependent and frequency-dependent transmission.

\noindent
\textit{(a) Exponential density-dependent transmission.} In this case, using \eqref{eq-km10} and \eqref{eq_kmA10b}, we obtain
\begin{equation}
\label{eq_kmA13}
\Lambda(t)=\overline{\beta} I(t)= \sum_{j=1}^{\infty}\overline{\beta}(1-\gamma)^{j-1}C(t-j).
\end{equation}

Comparing \eqref{eq_kmA13} and \eqref{eq-km7} we obtain
\begin{equation}
\label{eq_kmA14}
A_j=\overline{\beta}(1-\gamma)^{j-1},\qquad j=1,2,\ldots
\end{equation}

\noindent
\textit{(b) Exponential frequency-dependent transmission.} We now have

\begin{equation}
\label{eq_kmA15}
\Lambda(t)=\frac{\beta'}{N}I(t),
\end{equation}
so that
\begin{equation}
\label{eq_kmA16}
A_j=\frac{\beta'}{N}(1-\gamma)^{j-1},\qquad j=1,2,\ldots
\end{equation}
Let us interpret these results (compare with the discussion at the end of subsection \ref{subsec_sir_cum}). An infectious individual remains infectious from one step to the next with probability \(1-\gamma\) and is removed with probability $\gamma$. Thus, for an individual infected $j$ steps ago to contribute in the current step, they must not have been removed in the subsequent $j-1$ transitions. As we use the hypothesis of independence per step, the probability $\mathbb{P}$ of infectious survival is
\begin{equation}
\label{eq_kmA16b}
\mathbb{P}(\text{still in }I\text{ after }j-1\text{ steps})=(1-\gamma)^{j-1}.
\end{equation}
This is exactly the term that appears in (\eqref{eq_kmA14}) and (\eqref{eq_kmA16}): it is the survival function of the geometric distribution of infectious time (in steps), whose mean is \(1/\gamma\).

The discrete Kermack–McKendrick approach described above is useful for constructing generalizations of the SIR model that involve non-geometric removal. That is, generalizations that consist of abandoning the memoryless hypothesis (Bernoulli distribution with constant probability $\gamma$ per step). In terms of the discrete Kermack–McKendrick formalism, this translates directly into a change in the infectious survival factor that appears in the expression of $I(t)$ as a function of incidence $C(t)$, so that $A_k$ becomes proportional to a general survival function \cite{Diekmann2021}.

\section{Choice between Continuous and Discrete Modeling}
\label{sec_mod_cont_disc}

There are some fundamental  differences between continuous and discrete approaches. Here, we emphasize the discrete approach via difference equations, partly because it is less studied in the pedagogical literature. There are, however, many circumstances where the discrete approach is worth applying: 

\begin{itemize}
\item Epidemiological data are collected at discrete time intervals (days, weeks).
\item Day/night rhythms are often important in system dynamics \cite{Diekmann2025}.
\item We can model infection and removal (or recovery) as probabilities of events in the interval $(t,t+1]$. For example, in equation \eqref{eq-km2}, the term $\mathrm{e}^{-\Lambda(t)}$ is the probability of escape in that interval. For small values of $\Lambda$, the escape probability is $1-\Lambda(t)$.
\item Typical phenomena of dynamical systems (existence of stability regions and bifurcations) have a pedagogically clearer interpretation.
\end{itemize}

For more detailed comparisons between continuous and discrete approaches in compartmental epidemiological models, see, for example, references \cite{Allen1994, Allen2000, Wacker2020, Diekmann2021}.

\section{Python Implementations}
\label{sec_python}

\subsection{Basic discrete SIR model in Python}
\label{subsec_eq_dif_sir}

We will now numerically determine the discrete dynamics of the basic SIR model. The Python commands explained here can be found in \cite{Sasse_github}, in the file \texttt{sir\_basic\_discrete.ipynb}. We choose a model with frequency-dependent transmission, described by the system of equations \eqref{eq_discreta5a}, \eqref{eq_discreta5b}, and \eqref{eq_discreta5c}. An implementation of the basic discrete SIR model with density-dependent transmission can be obtained by making minor modifications to the Python commands described below.

The system of difference equations \eqref{eq_discreta5a}, \eqref{eq_discreta5b}, and \eqref{eq_discreta5c}, although theoretically complicated because it is nonlinear, is relatively simple to implement from a numerical point of view, as the quantities $S_{n+1}$, $I_{n+1}$, and $R_{n+1}$ are completely determined in terms of $S_{n}$, $I_{n}$, and $R_{n}$, as well as the parameters $\beta$ and $\gamma$. We can solve such a system iteratively in Python, as we will show below. The commands described can be reproduced in \textit{Google Colab} or \textit{Jupyter Notebook}. We use only elementary programming resources, without concern for algorithm optimization.

Initially, we load the libraries that will be used:

\begin{tcolorbox}
\begin{verbatim}
import numpy as np
import matplotlib.pyplot as plt
\end{verbatim}
\end{tcolorbox}

We will use, as an example, a typical case of measles spread in a closed, non-immune population \cite{Anderson1991}. In this case, $\beta'=1.5\,\mbox{days}^{-1}$ and $\gamma=0.12\,\mbox{days}^{-1}$, so that $\mathcal{R}_0^{\prime}=\beta'/\gamma=12.5$. These constants can be changed later based on observations. Furthermore, we  take initial values $S=4999$, $I=1$, and $R=0$, so that $N=5000$.

In the next cell, we define the parameters $\beta'$ and $\gamma$ and the total simulation length:
\begin{tcolorbox}
\begin{verbatim}

# Model parameters

beta = 1.5          # Infectiousness
gamma = 0.12        # Recovery rate
T = 50              # Simulation length (days)
\end{verbatim}
\end{tcolorbox}
\noindent
The initial conditions are defined below:
\begin{tcolorbox}
\begin{verbatim}

# Initial conditions

S0 = 4999           # Initial number of susceptibles
I0 = 1              # Initial number of infectious
R0 = 0.0            # Initial number of recovered
\end{verbatim}
\end{tcolorbox}
\noindent
Although not mandatory, it is good computational practice to predefine the sizes of the vectors (1D arrays) to be used, if possible. That is, we define each vector with zero components in a quantity equal to the number of points to be calculated, including the initial point:
\begin{tcolorbox}
\begin{verbatim}

#Variable initialization

steps = T + 1       # Includes the initial point
S = np.zeros(steps)
I = np.zeros(steps)
R = np.zeros(steps)
time = np.linspace(0, T, steps)  # Time vector
\end{verbatim}
\end{tcolorbox}

The initial conditions correspond to the values of the zero components of the 1-dimensional arrays:
\begin{tcolorbox}
\begin{verbatim}

#Initial state

S[0] = S0
I[0] = I0
R[0] = R0
N = S0 + I0 + R0
\end{verbatim}
\end{tcolorbox}
\noindent

We can now perform the iterative calculations:
\begin{tcolorbox}
\begin{verbatim}

#Iteration with difference equations

for t in range(steps - 1):
  S[t+1] = S[t] - beta * S[t] * I[t] / N
  I[t+1] = I[t] + (beta * S[t] * I[t] / N - gamma * I[t])
  R[t+1] = R[t] + gamma * I[t]
\end{verbatim}
\end{tcolorbox}

To visualize the results, simply execute the commands:
\begin{tcolorbox}
\begin{verbatim}
plt.figure(figsize=(7, 4), dpi=200)
plt.plot(time, S/N, label="Susceptibles ($S$)", lw=1.5)
plt.plot(time, I/N, label="infectious ($I$)", lw=1.5)
plt.plot(time, R/N, label="Recovered ($R$)", lw=1.5)
plt.title("Basic discrete SIR model")
plt.xlabel("$t$ (days)")
plt.ylabel("Proportion of population")
plt.legend()
plt.grid()
plt.savefig("sir_basic_discrete.png", dpi=300) # Save in high resolution
plt.show()
\end{verbatim}
\end{tcolorbox}

\begin{figure}[ht]
\centering
\includegraphics[width=0.8\linewidth]{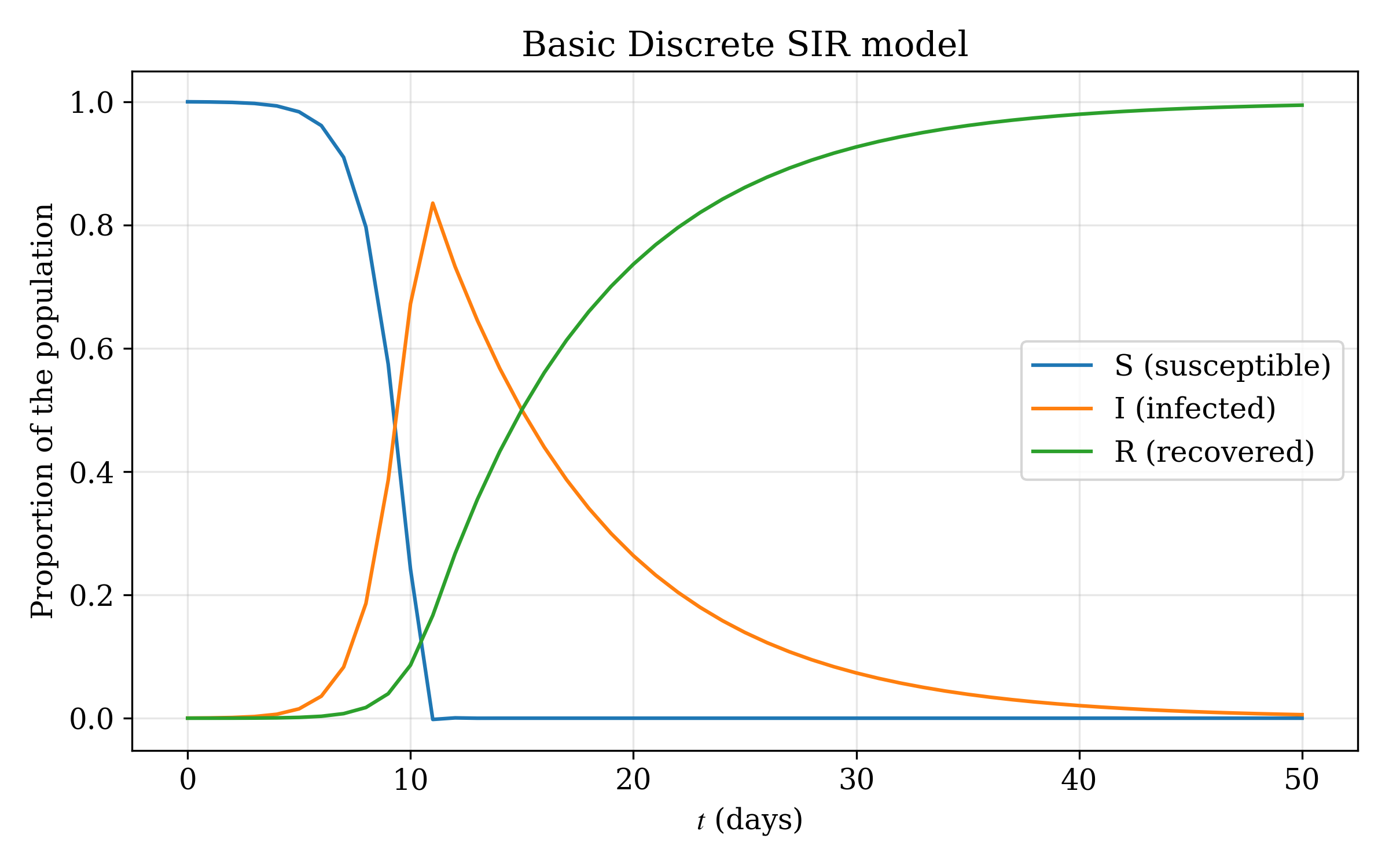}
\caption{Basic discrete SIR model: simulation based on difference equations.}
\label{sir_basic_discrete}
\end{figure}

The result is shown in Figure \ref{sir_basic_discrete}.

\FloatBarrier
We can note that the maximum number of infectious individuals reaches approximately 83\% of the population on the eleventh day. The epidemic is nearly over after approximately 50 days.

\subsection{Basic continuous SIR model in Python}
\label{subsec_eq_diferenciais_sir}

We will now construct the basic continuous SIR model with frequency-dependent transmission. The Python commands explained here are available in \cite{Sasse_github}, in the file \\ \texttt{sir\_basic\_cont.ipynb}. 

The model is now represented by the set of nonlinear differential equations defined by \eqref{eq-continua1b}, \eqref{eq-continua2b}, and \eqref{eq-continua3b}. To compare results, we will use the same parameters and initial conditions as in the discrete model. Such a system of equations can be solved numerically in Python, as we will show below.

Initially, we load the necessary libraries:

\begin{tcolorbox}
\begin{verbatim}
import numpy as np
import matplotlib.pyplot as plt
from scipy.integrate import solve_ivp
\end{verbatim}
\end{tcolorbox}

We define the model parameters:

\begin{tcolorbox}
\begin{verbatim}

Model parameters

beta = 1.5           # Infectiousness
gamma = 0.12         # Removal rate
T = 50               # Total duration (in days)
\end{verbatim}
\end{tcolorbox}
\noindent
and the initial conditions:

\begin{tcolorbox}
\begin{verbatim}
S0 = 4999            # Initial number of susceptibles
I0 = 1               # Initial number of infectious
R0 = 0.0             # Initial number of recovered
N = S0 + I0 + R0     # Total number of individuals
IC = [S0, I0, R0]    # List of initial conditions
\end{verbatim}
\end{tcolorbox}

We define the system of equations for the model:

\begin{tcolorbox}
\begin{verbatim}

System of equations

def sir(t, y):
    S, I, R = y
    dSdt = -beta * S * I / N
    dIdt = beta * S * I / N - gamma * I
    dRdt = gamma * I
    return [dSdt, dIdt, dRdt]
\end{verbatim}
\end{tcolorbox}

Let us choose a grid with 500 points for the time interval $[0,50]$:

\begin{tcolorbox}
\begin{verbatim}
time = np.linspace(0, T, 500)
\end{verbatim}
\end{tcolorbox}

To solve the system numerically, we will use the Runge-Kutta method:

\begin{tcolorbox}
\begin{verbatim}
solution = solve_ivp(sir, [0, T], IC, t_eval=time, method="RK45")
S, I, R = solution.y
\end{verbatim}
\end{tcolorbox}

\begin{figure}[ht]
\centering
\includegraphics[width=0.8\linewidth]{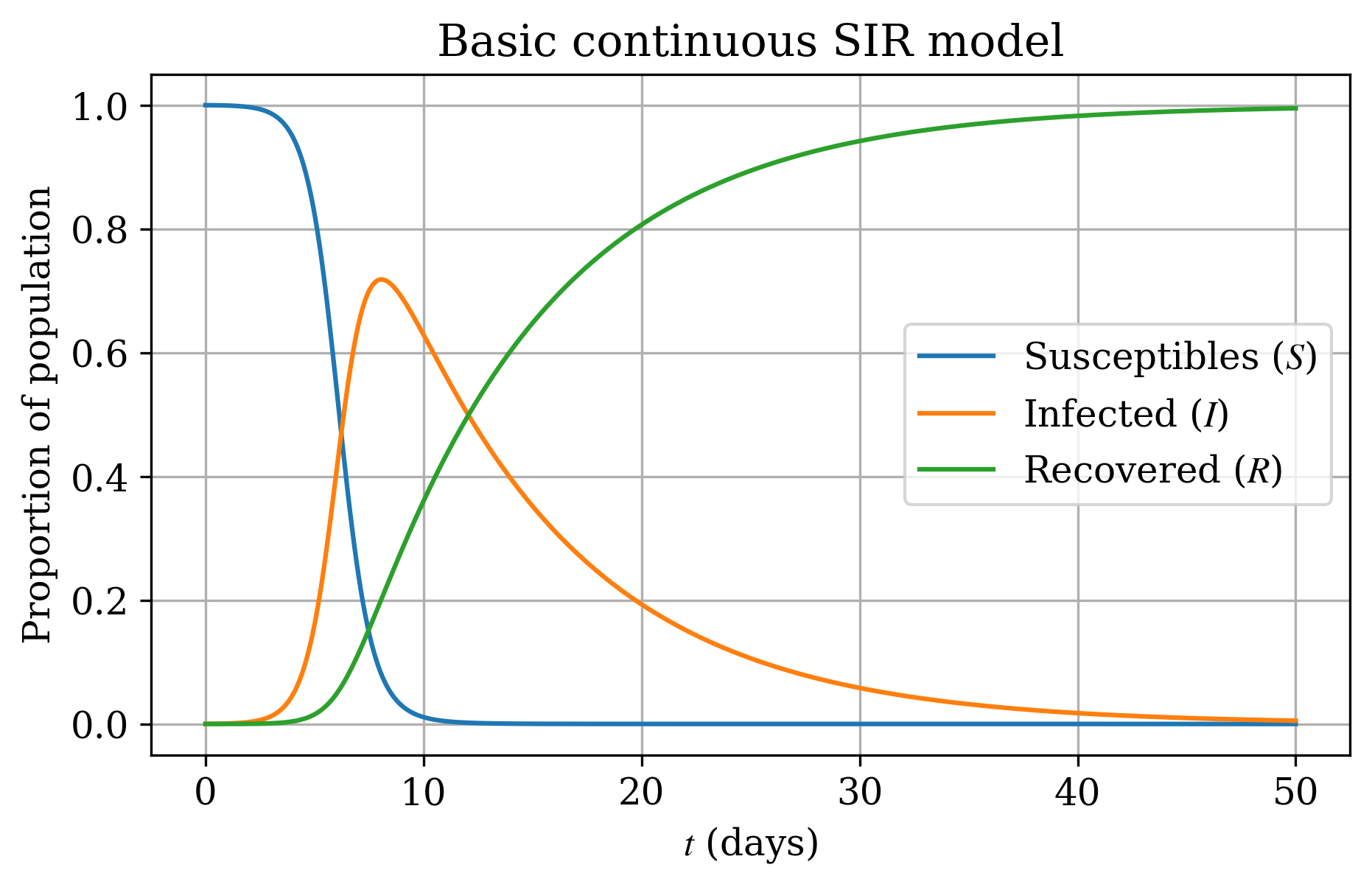}
\caption{Basic continuous SIR model: simulation based on differential equations.}
\label{sir_basic_cont}
\end{figure}

We can now graphically visualize the results:

\begin{tcolorbox}
\begin{verbatim}
plt.figure(figsize=(7, 4), dpi=200)
plt.plot(time, S/N, label= "Susceptibles (S)", lw=1.5)
plt.plot(time, I/N, label= "infectious (I)", lw=1.5)
plt.plot(time, R/N, label= "Recovered", lw=1.5)
plt.title("Basic continuous SIR Model")
plt.xlabel("$t$ (days)")
plt.ylabel("Proportion of population")
plt.legend()
plt.grid()
plt.savefig("sir_cont1.png", dpi=300) # Save in high resolution
plt.show()
\end{verbatim}
\end{tcolorbox}
\noindent
The result is presented in Figure \ref{sir_basic_cont}.

\FloatBarrier

This graph is qualitatively similar to that obtained with difference equations, as shown in Figure \ref{sir_basic_discrete}. However, the continuous model predicts that the maximum number of infectious individuals is approximately 72\% of the population on the eighth day. As in the discrete model, here the epidemic is nearly over after approximately 50 days. The graphs suggest that this model has a disease-free equilibrium state. This can be verified by setting $dS/dt=dI/dt=dR/dt = 0$ in \eqref{eq-continua_freq1}-\eqref{eq-continua_freq3} as discussed in \cite{Blackwood2018}.

\subsection{SIR model with exponential transmission in Python}

\begin{figure}[ht]
\centering
\includegraphics[width=0.8\linewidth]{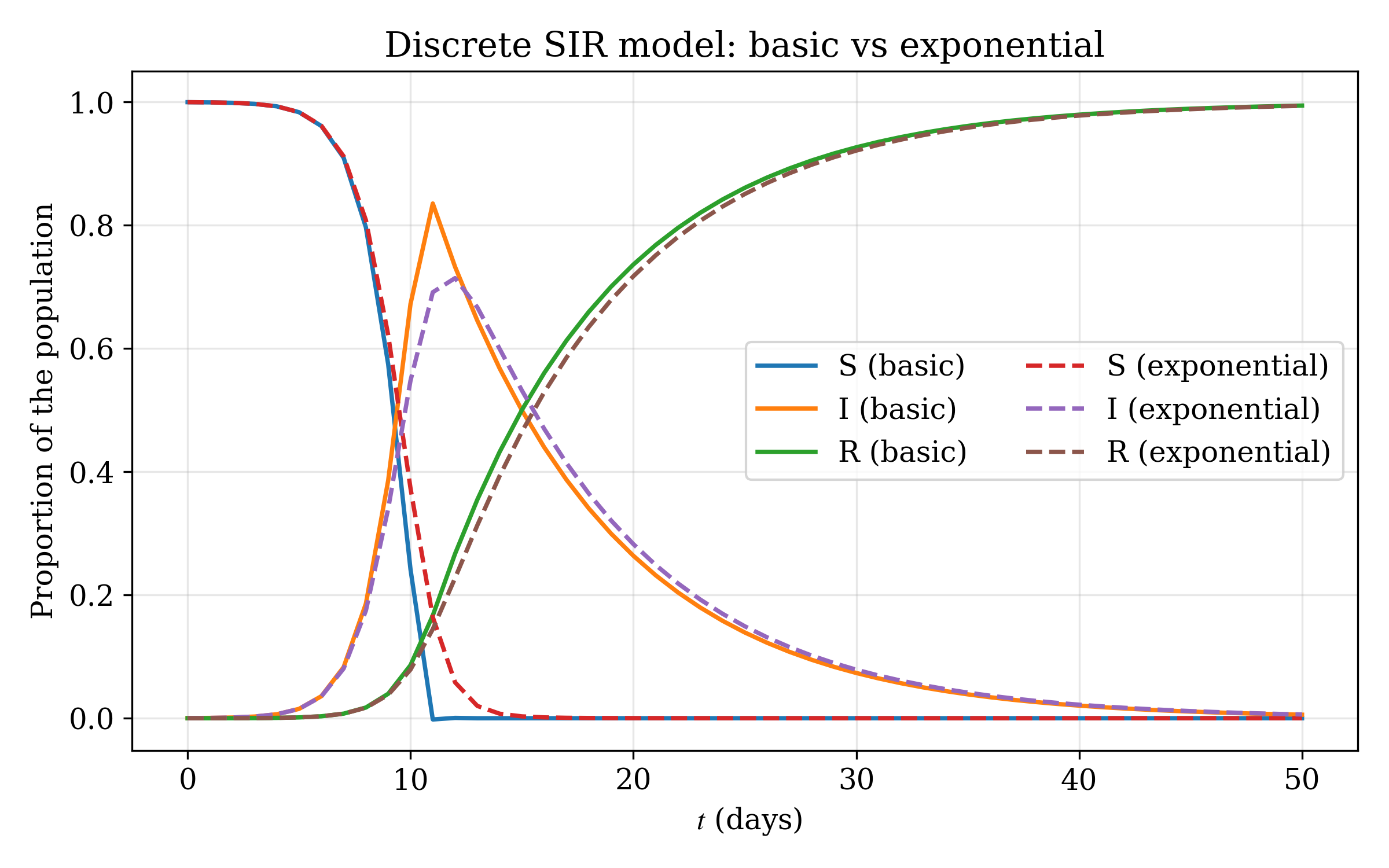}
\caption{Basic and exponential discrete SIR models.}
\label{fig_sir_exp}
\end{figure}

The SIR model with exponential transmission, defined by \eqref{eq-km23b}, \eqref{eq-km24b}, and \eqref{eq-km25b}, can be found in the repository \cite{Sasse_github}, file \texttt{sir\_exp\_disc.ipynb}. The results, showing the graphs for the exponential and basic SIR models, are presented in Figure \ref{fig_sir_exp}.
\FloatBarrier

We note that the exponential model ``smooths'' the graphs generated by the basic model. We should note that the basic SIR model can generate negative values for $S_i$.

\section{Parameter Estimation}
\label{sec_estimation}

In this section, we present a method for estimating the parameters of our models.  We  take as an example the basic discrete SIR model with frequency-dependent transmission, according to \eqref{eq_discreta5a}–\eqref{eq_discreta5c}. Let us assume that we have observational data in the form of a time series:
\begin{equation}
I_{\mathrm{obs}}(0), I_{\mathrm{obs}}(1), \ldots, I_{\mathrm{obs}}(n),
\end{equation}
with the number of active infectious individuals $I(t)$, for each discrete time unit, from 0 to $n$. It is important to distinguish \textit{prevalence} (active infectious) from \textit{incidence} (new cases in the interval) and \textit{cumulative cases}. The procedure described here assumes $I_{\mathrm{obs}}(n)$ as prevalence data. If the available data are daily incidence, it is necessary to adapt the method to fit incidence, not prevalence. From this, we define the initial conditions of the model:
\begin{equation}
I_0 := I_{\mathrm{obs}}(0),
\qquad
S_0 := N - I_0 - R_0.
\end{equation}
If we have a pair of values $(\beta',\gamma)$ given, we can use \eqref{eq_discreta5a}–\eqref{eq_discreta5c} to generate a predicted trajectory
\begin{equation}
\Big(S_{\mathrm{pred}}(k), I_{\mathrm{pred}}(k), R_{\mathrm{pred}}(k)\Big),
\quad k=0,1,\ldots,n.
\end{equation}
A simple method to estimate the optimal value of $(\beta', \gamma)$,  from $I$ data consists of minimizing the sum of squared residuals (least-squares method):
\begin{equation}
\label{eq_obj_ls}
J(\beta', \gamma)
=
\sum_{k=0}^{n} \left[I_{\mathrm{pred}}(k;\beta',\gamma)-I_{\mathrm{obs}}(k)\right]^2.
\end{equation}
The estimation problem becomes:
\begin{equation}
\label{eq_estimacao_min}
(\widehat{\beta'}, \widehat{\gamma})
=
\arg\min_{\beta'\ge 0,\ 0\leq\gamma\leq 1} J(\beta',\gamma).
\end{equation}

The fitting process described above, based only on $I(t)$, is useful in the pedagogical context and can work well in series with  a clearly observable peak and decay. However, in real data, estimation can be sensitive to (i) underreporting, (ii) reporting delays, (iii) temporal changes in behavior and control measures (which make $\alpha$ effectively variable over time), and (iv) uncertainties in initial conditions, particularly $R_0$ and even $I_{\mathrm{obs}}(0)$. In more realistic inferential applications, it may be necessary to use observation models (e.g., Poisson or negative binomial) and maximum likelihood or Bayesian methods \cite{Keeling2008, Vynnycky2010}. 

As an illustrative dataset, we use daily observations from the 1978 influenza outbreak in an English boarding school, involving a closed population of 763 students \cite{Anonymous1978}. The data report the daily number of students confined to bed during the outbreak. As a simplifying assumption, we assume that the daily number of students confined to bed works as a good proxy for the number of active infectious individuals.

\begin{figure}[ht]
\centering
\includegraphics[width=0.8\linewidth]{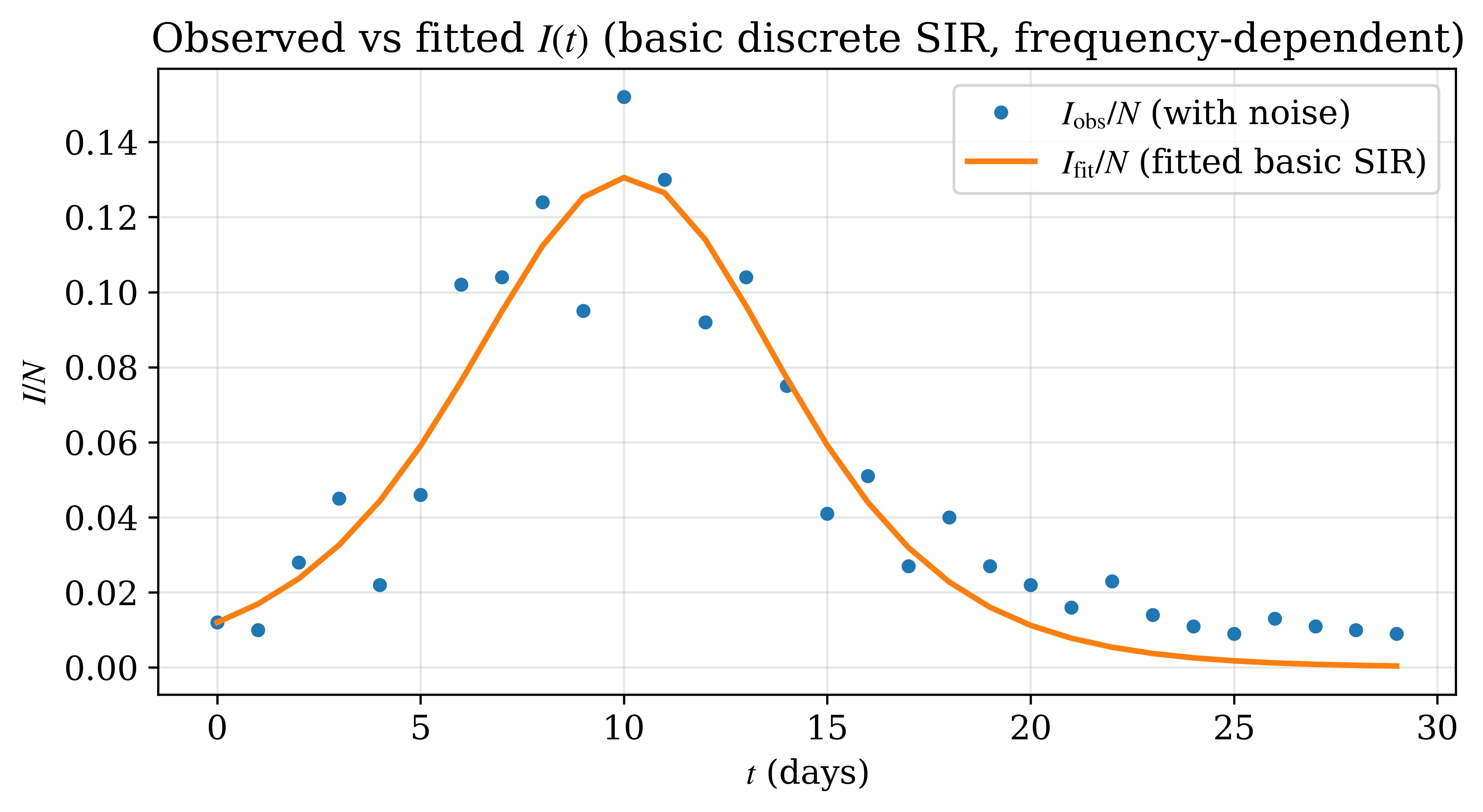}
\caption{Fitting the basic discrete SIR model to data to a school influenza outbreak.}
\label{fig-sir_basic_fit_I}
\end{figure}

A Python implementation of the method described here for the basic discrete SIR model with frequency-dependent transmission can be found in \cite{Sasse_github}, in the file \linebreak \texttt{sir\_basic\_disc\_est.ipynb}.  The resulting observed-versus-fitted plot is shown in Figure \ref{fig-sir_basic_fit_I}. 
The plot of the residuals is in the output file. We note that the fit tracks the peak well but sits systematically slightly below the data from about day 18 onward.
\FloatBarrier

\section{Conclusions}
\label{sec_conclusoes}

In this article, we present a didactic approach to compartmental models, with special emphasis on the basic hypotheses and the fundamental aspects of transport between compartments, which are often taken for granted in the existing educational literature. We focus on the SIR model, which, despite being one of the simplest, is sufficient for our analysis. The discussion can be generalized to more complicated models.

We stress that most  pedagogical works on compartmental models do not discuss the construction of transmission mechanisms in detail, nor do they discuss the differences between density-dependent and frequency-dependent transmission. In our presentation, such a discussion plays a central role.

We argue that the discrete approach should not be seen only as a numerical discretization of the continuous model but also as a natural form of modeling when data are observed at discrete times and when one wishes to interpret infection and recovery as probabilistic events. Additionally, we present a more general version of the model, with an exponential term associated with the cumulative force of infection, which has a probabilistic interpretation of escape probability and tends to keep variables in coherent domains while conceptually connecting to the Kermack–McKendrick formalism \cite{Breda2012, Diekmann2021}. In a pedagogical context, this comparison provides an opportunity to discuss the difference between \textit{approximation} and \textit{model}, as well as the importance of making hypotheses and units of measurement explicit.

The Python implementations presented in the article are not optimized and are intended to promote accessibility. In the classroom, they allow students to explore parameters, identify thresholds associated with the basic reproduction number, and build intuition about the role of transmission mechanisms. The proposed activities can be used in courses that integrate computational experimentation. This article can serve as a basis for investigating more complex compartmental models, involving different compartments and modifications of the SIR model, such as the introduction of a reporting parameter to estimate the proportion of underreported cases \cite{Kalachev2024}, nonlinear incidence \cite{Habott2024}, time-dependent parameters \cite{Chen2020}, and age structuring \cite{Ram2021}.

Future work involves the inclusion of new compartmental groups and heterogeneities (age, spatial, and social structure), the generalization of the geometric kernel $A_k$ (which results in the exponential SIR models) to more general forms, the systematic comparison between deterministic and stochastic models, and the discussion of parameter estimation methods from real data, with special attention to the model-dependent nature of quantities such as $\mathcal{R}_0$ \cite{Li2011, Delamater2019, Lim2020}.

\section*{Conflict of Interest}
The authors declare that they have no known competing financial interests or personal relationships that could have appeared to influence the work reported in this paper.

\section*{Acknowledgments}
The authors gratefully acknowledge the support of the Department of Mathematics at the State University of Santa Catarina (UDESC).

\end{document}